\documentstyle[12pt]{article}
\def\tr{{\rm T}\!{\rm r}}

\def\CA{{\cal A}}    
 \def\CF{{\cal F}}  \def\CH{{\cal H}} 
   \def\CL{{\cal L}} 
 \def\CN{{\cal N}}  \def\CP{{\cal P}}

   \def\BD{{\bf D}} 
    
 \def\BJ{{\bf J}}  \def\BL{{\bf L}} 
 \def\BN{{\bf N}}  
  \def\BS{{\bf S}} \def\BT{{\bf T}}

\def\Ba{{\bf a}} \def\Bb{{\bf b}}

  \def\Bo{{\bf o}} 
   
   \def\Bx{{\bf x}}

\newfont{\frak}{eufm10 scaled 1200}

\def\half{{1\over2}}

\def\complex{{\mathchoice
{\setbox0=\hbox{$\displaystyle\rm C$}\hbox{\hbox to0pt
{\kern0.4\wd0\vrule height0.9\ht0\hss}\box0}}
{\setbox0=\hbox{$\textstyle\rm C$}\hbox{\hbox to0pt
{\kern0.4\wd0\vrule height0.9\ht0\hss}\box0}}
{\setbox0=\hbox{$\scriptstyle\rm C$}\hbox{\hbox to0pt
{\kern0.4\wd0\vrule height0.9\ht0\hss}\box0}}
{\setbox0=\hbox{$\scriptscriptstyle\rm C$}\hbox{\hbox to0pt
{\kern0.4\wd0\vrule height0.9\ht0\hss}\box0}}}}
\def\Co{{\mathchoice
{\setbox0=\hbox{$\displaystyle\rm C$}\hbox{\hbox to0pt
{\kern0.4\wd0\vrule height0.9\ht0\hss}\box0}}
{\setbox0=\hbox{$\textstyle\rm C$}\hbox{\hbox to0pt
{\kern0.4\wd0\vrule height0.9\ht0\hss}\box0}}
{\setbox0=\hbox{$\scriptstyle\rm C$}\hbox{\hbox to0pt
{\kern0.4\wd0\vrule height0.9\ht0\hss}\box0}}
{\setbox0=\hbox{$\scriptscriptstyle\rm C$}\hbox{\hbox to0pt
{\kern0.4\wd0\vrule height0.9\ht0\hss}\box0}}}}
\def\Rl{{\mathchoice
{\setbox0=\hbox{$\displaystyle\rm R$}\hbox{\hbox to0pt
{\kern0.4\wd0\vrule height0.9\ht0\hss}\box0}}
{\setbox0=\hbox{$\textstyle\rm R$}\hbox{\hbox to0pt
{\kern0.4\wd0\vrule height0.9\ht0\hss}\box0}}
{\setbox0=\hbox{$\scriptstyle\rm R$}\hbox{\hbox to0pt
{\kern0.4\wd0\vrule height0.9\ht0\hss}\box0}}
{\setbox0=\hbox{$\scriptscriptstyle\rm R$}\hbox{\hbox to0pt
{\kern0.4\wd0\vrule height0.9\ht0\hss}\box0}}}}

\def\be{\begin{equation}}
\def\ee{\end{equation}}
\def\bea{\begin{eqnarray}}
\def\eea{\end{eqnarray}}

\begin{document} 

\begin{titlepage}
\vspace*{-4ex}

\null \hfill Preprint TU-525  \\
\null \hfill October 1997 \\[4ex]

\begin{center}

\bigskip
\bigskip

{\Large \bf
Differential Calculus on Fuzzy Sphere and Scalar Field
}\\[5em]

Ursula Carow-Watamura\footnote{
e-mail: ursula@tuhep.phys.tohoku.ac.jp} \  and \  
Satoshi Watamura\footnote{e-mail: watamura@tuhep.phys.tohoku.ac.jp}

Department of Physics \\
Graduate School of Science \\
Tohoku University \\
Aoba-ku, Sendai 980-77, JAPAN \\ [2ex]
\end{center}

\bigskip

\begin{abstract}
We find that there is an alternative possibility to define the 
chirality operator on the fuzzy sphere, 
due to the ambiguity of the operator ordering.
Adopting this new chirality operator and the corresponding Dirac operator, 
we define Connes' spectral triple on the fuzzy sphere 
and the differential calculus. The differential calculus based on this new 
spectral triple is simplified considerably. 
Using this formulation the action of
the scalar field is derived.
\end{abstract}

\end{titlepage}

\eject  


\section{Introduction}
Recent developments of the nonperturbative aspects in string theory 
give signs that the noncommutative geometry is 
the right way to describe the quantum regime of spacetime \cite{BFSS,IKKT}.
Compared to other attempts in the past to generalize the 
spacetime structure, the noncommutative geometry has the advantage that it 
can preserve the kinematical properties of the geometry treated, and 
the framework is general enough to extend to arbitrary dimension and
 to include nontrivial topology.

The notion of 'noncommutative geometry' has been first considered by 
von Neumann. Later Connes extended this notion to the 
'noncommutative differential geometry' which means that we consider 
the noncommutative function algebra associated with a certain 
geometry together with a differential calculus. 
Technically, it is the translation of the
various geometrical concepts into the algebraic terminology, where the 
original space $X$ is replaced by a function algebra 
$C^\infty(X)$ of the smooth functions on $X$.
Once the geometry described by purely algebraic language, it is possible 
to replace its 
structure algebra by a noncommutative algebra 
and we obtain the geometry of the noncommutative space \cite{ConnesNCG} . 

The fuzzy sphere is one of the simplest examples of noncommutative 
geometry, and has been discussed by many authors. 
The algebra of the fuzzy sphere appeared already in many 
different contexts \cite{Haldane83,FanoOrto86,Hoppe89,DeWHoNi88,MadoreCQG92}.  
This algebra can be constructed 
by applying Berezin's quantization \cite{Berezin74,Berezin75a} 
to the function algebra over the 
sphere, and the result is a matrix algebra of finite 
dimension. 
However the way to introduce the differential 
calculus is not unique even if we require that in the commutative 
limit the standard differential calculus is reproduced.
Proposals for differential calculi have been made by 
\cite{ConnesNCG,DuboKern90} and especially for the 
quantum sphere in 
\cite{MadoreCQG92,GrosMado92,GrosKlim96}. 

In this paper we propose a new spectral triple for the differential 
calculus on the fuzzy sphere 
within Connes' framework. As is known, in order to apply 
Connes' construction of the differential algebra, it is necessary 
to define the Dirac operator and the chirality operator
\cite{GrosPres95,GrosKlim96a,CW96}. 
The approach taken here to construct the Dirac operator 
is similar to ref.\cite{CW96}. The difference is that we adopt 
a different chirality operator, which is possible due to the 
ambiguity of operator ordering.
Using this new chirality operator and the corresponding Dirac operator, 
we define a differential algebra 
within Connes' framework which has a rather simple structure. 
With the resulting 
differential calculus we derive 
the action of the scalar boson field on the noncommutative sphere. 

This paper is organized as follows. In section 2, we give some definitions
which are necessary to define the chirality operator and the Dirac operator
on the fuzzy sphere. In section 3, the new chirality operator and 
Dirac operator are introduced and the differential calculus is
defined using these operators. In section 4, we apply this 
differential calculus to the scalar field and derive the action of the 
scalar boson.
Section 5 is the discussion.

\section{Algebra of Fuzzy Sphere}

Let us first recall the definition of the fuzzy sphere. 
The algebra $\CA_N$ 
of the fuzzy sphere is generated by the operators $\Bx_i$ ($i=1,2,3$) 
satisfying 
\be
[\Bx_i,\Bx_j]=i\alpha \epsilon_{ijk}\Bx_k \ ,
\ee
with the constraint
\be
\Bx_i\Bx_i=\ell^2 \ .\label{radius}
\ee
Each operator of $\CA_N$ is 
represented by a matrix acting on the $N+1$-dimensional
Hilbert space $\CF_N$. The "Planck constant" $\alpha$ is 
a central element,
$[\alpha , \Bx_i]=0$, and its value is determined by eq.(\ref{radius}) 
as $\alpha = {2\ell\over \sqrt{N(N+2)}}$. $\alpha$ vanishes 
in the limit $N\rightarrow \infty$, 
which is the commutative limit of the algebra.
Apparently, the matrices $\Bx_i$ can be identified
 with the generators of the $su(2)$ Lie algebra
and the generated algebra is equivalent to the algebra of 
$(N\!+\!1)\!\times\!(N\!+\!1)$
matrices, 
$M_{\complex}(N+1)$. 

A natural object as a field on the fuzzy sphere is an $\CA_N$-bimodule, where
 we can consider left multiplication and right multiplication of 
the algebra $\CA_N$ 
onto this field.  In the commutative limit, the distinguishing of right and 
left multiplication is not important, however in the noncommutative case
it has rather nontrivial consequences.
Since left multiplication and right multiplication commute even 
in the case of  noncommutative space, 
the $\CA_N$-bimodule can be considered as a left module over the algebra
 $\CA_N\otimes \CA_N^{\Bo}$, 
where $\CA_N^{\Bo}$ denotes the opposite algebra which is defined by:
\be
\Bx_i^{\Bo}\Bx_j^{\Bo}\equiv (\Bx_j\Bx_i)^{\Bo} \ , \ 
\Bx_i\in \CA_N \ .
\ee
The action of $\Ba ,\Bb \in\CA_N$ 
onto the $\CA_N$-bimodule
$\Psi\in\Gamma\!\CA_N$ is
\be
\Ba\Bb^{\Bo}\,\Psi\equiv \Ba\,\Psi\,\Bb \ .
\ee
Therefore, when we consider operators acting on fields 
over the fuzzy sphere,
 these operators are elements of the algebra $\CA_N\otimes \CA_N^{\Bo}$
\cite{DuboMass95}. 
And it is natural to define a differential operator acting on the field
 as an element in $\CA_N\otimes \CA_N^{\Bo}$ .
As we shall see in the following,
 the Dirac operator derived in ref.\cite{CW96} also fits to this scheme.

\section{Dirac Operator and Differential Calculus}

Introducing the spinor field $\Psi$ as an $\CA_N$-bimodule
$\Gamma\CA_N\equiv\complex^2\otimes \CA_N$
we define the chirality operator 
and the Dirac operator $\BD$ as operators on the spinor field,
i.e. as $2\times2$ matrices the entries of which are 
elements of the algebra $\CA_N\otimes \CA_N^{\Bo}$.
The construction of the Dirac operator can be performed 
by the following steps:

\begin{itemize}
\item[(a)] define a chirality operator which has a standard commutative limit,\item[(b)] define the Dirac operator requiring that 
it anticommutes with the chirality operator and, 
 in the commutative limit it reproduces the standard Dirac
operator on the sphere.
\end{itemize}

In ref.\cite{CW96} we have constructed the Dirac operator by applying 
this approach. There, we have introduced a chirality operator given by
\be
\gamma_\chi={1\over \CN}(\sigma_i\otimes \Bx_i+{\alpha\over 2})\ ,
\ee
where $\CN$ is a normalization constant determined by the condition
$\gamma_\chi^2=1$, and
$\sigma_i$ ($i=1,2,3$) are the Pauli matrices.
In the commutative limit, the operator $\Bx_i$ can be identified with
the homogeneous coordinate $x_i$ of sphere and 
the chirality operator given in eq.(5) becomes ${\sigma_ix_i\over\ell}$
 which is the standard
chirality operator invariant under the rotation. We have also defined 
the corresponding Dirac operator. 
Since the chirality operator does not commute
with the original algebra $\CA_N$,
we proposed a modified algebra which commutes with $\gamma_\chi$,
to obtain Connes' spectral triple. 

In this paper,  we use a different chirality operator
which also produces the standard chirality operator in 
the commutative limit. This is possible due to the fact that
the condition for the commutative limit of 
the chirality operator does not fix its noncommutative generalization 
uniquely, i.e. we have an ambiguity in 
the operator ordering.
The new chirality operator commutes with the original algebra $\CA_N$
and we obtain the spectral triple with the original 
algebra $\CA_N$; hence, a modification of the algebra as 
proposed in ref.\cite{CW96} is not necessary. 

Requiring the same condition as for $\gamma_\chi$, i.e. (a),
the new chirality operator is determined as
\be
\gamma_\chi^{\Bo}={1\over \CN}(\sigma_i\otimes\Bx^{\Bo}_i-{\alpha\over 2}) \ .
\ee
In the following we do not explicitly write the symbol $\otimes$ since its
meaning is obvious.
The normalization constant $\CN$ is defined 
by the condition, 
\be
(\gamma_\chi^{\Bo})^2=1 \ ,
\ee
as $\CN={\alpha \over 2}(\BN+1)$.

Searching the Dirac operator $\BD$ by requiring 
the condition (b), i.e., $\{\gamma_\chi^{\Bo},\BD\}=0$, we find
\be
\BD={i\over \ell\alpha}
\gamma_\chi^{\Bo}\epsilon_{ijk}\sigma_i\Bx_j^{\Bo}\Bx_k\ .
\label{Dirac0}
\ee
Note that this Dirac operator is selfadjoint,  
$\BD^\dagger = \BD$.

Acting with this operator on the spinor $\Psi\in\Gamma\CA_N$,
we obtain
\be
\BD\Psi={i\over \ell}\gamma_\chi^{\Bo}\chi_i\BJ_i\Psi\ ,
\label{Dirac1}
\ee
where
\be
\BJ_i=\BL_i+\half\sigma_i \ ,
\ee
and
\be
\chi_i\equiv\epsilon_{ijk}\Bx_j\sigma_k
\ee
The angular momentum operator is defined by the adjoint action
of $\Bx_i$:
\be
\BL_i\Psi\equiv {1\over\alpha}[\Bx_i,\Psi]={1\over\alpha}(\Bx_i\Psi-\Psi\Bx_i)
\ee

If we replace the chirality operator 
$\gamma_\chi^{\Bo}$ by $\gamma_\chi$ in 
eq.(\ref{Dirac1}), we obtain the Dirac operator 
 given in ref.\cite{CW96}, since $\chi_i\BJ_i=\Sigma-\chi$, 
where $\Sigma=-i\epsilon_{ijk}\sigma_i\Bx_j\BL_k$, $\chi=\sigma_i\Bx_i$, 
and $\gamma_\chi(\Sigma-\chi)$ is 
the Dirac operator given in ref.\cite{CW96} 
up to a normalization constant. 
The essential part of the Dirac operator 
in eq.({\ref{Dirac0}) is the factor 
$\epsilon_{ijk}\sigma_i\Bx_j\Bx^{\Bo}_k$, which anticommutes with 
both  $\gamma_\chi^{\Bo}$ and $\gamma_\chi$.

The second condition of (b) concerning the commutative limit of the Dirac 
operator is also satisfied. 
If we take the commutative limit of each operator 
$\chi_i, \BJ_i, \gamma^{\Bo}_\chi$ in eq.(\ref{Dirac1}), we obtain

\be
\BD_\infty={i\over\ell}\gamma_\chi\chi_i\BJ_i
={i\over \ell^2}(\sigma_l x_l)\epsilon_{ijk} x_i\sigma_j (iK_k+\half\sigma_k)
=-{1\over\ell} (i\sigma_i K_i+1)
\ee
where $x_i$ is the homogeneous coordinate of $S^2$ and $K_i$ is the 
Killing vector. Therefore, in the commutative limit 
this Dirac operator is equivalent to the standard Dirac operator.
For details see also ref.\cite{CW97a}.

In order to establish Connes' triple we have 
to identify the Hilbert space. 

The space of the fermions $\Psi\in\CA_n\otimes M_2(\complex)$ 
defines the Hilbert space $\CH_N$ with the norm 
\be
<\Psi|\Psi>=\tr_{\CF}(\Psi^\dagger\Psi)
\ee 
where $\tr_{\CF}$ is the trace over the $(N+1)$ dimensional 
Hilbert space
$\CF_N$.

The dimension of the Hilbert space $\CH_N$ is $2(N+1)^2$ and the 
trace over
$\CH_N$ is the trace 
over the spin suffices and over the $(N+1)^2$ 
dimensional space of the matrices.
Since the Dirac operator is defined in the algebra 
$\CA_N\otimes\CA_N^{\Bo}$,
in general the trace must be taken 
for operators of the form $\Ba\Bb^{\Bo}$, 
with $\Ba,\Bb\in\CA_N$, and it is given by:
\be
\tr_{\CH}\{\Ba\Bb^{\Bo}\}=\sum^{2(N+1)^2}_{J=1}
<\Psi_J|\Ba\Bb^{\Bo}|\Psi_J>
=2\tr_{\CF}\{\Ba\}\tr_{\CF}\{\Bb\} \ ,\label{trace}
\ee
where $\Psi_J$ is an appropriate basis in $\CH_N$. The factor $2$ 
on the r.h.s. comes from the sum over spin suffices.

To examine the structure of the Hilbert space we compute the 
spectrum $\lambda_j$ of the Dirac operator:
\be
\BD^2\Psi_{jm}=\lambda_j^2\Psi_{jm}
\ee
where $\Psi_{jm}$ is a state with total angular momentum $j$, 
$\BJ^2\Psi_{jm}=j(j+1)\Psi_{jm}$, and $\BJ_3\Psi_{jm}=m\Psi_{jm}$.
$j$ and $m$ are half integer and run $\half\leq j\leq N+\half$,
$-j\leq m\leq j$.
The spectrum is given by 
\be
\lambda^2_j=(j+\half)^2[1
+{1-(j+\half)^2\over N(N+2)}]\ .
\ee
This spectrum is equivalent to the spectrum of the Dirac operator 
given in ref.\cite{CW96}. 
It corresponds to the classical spectrum in the limit $N\rightarrow \infty$, 
except for the part with maximal angular momentum $j=N+\half$.  
When the angular momentum takes its maximal value 
we see that $\lambda_{N+\half}=0$.  This happens since there is no
chiral pair for the spin $N+\half$ state and therefore this part must be 
a zeromode for consistency. We can also confirm this property by 
computing the index, $\tr_{\CH}(\gamma_\chi)$. 
Since these zeromodes have no classical analogue, one way 
to treat this situation is to project them out from the Hilbert space. 
On the other hand, as we shall see, 
the contribution of the zeromodes in the integration 
is of order ${1\over N}$ and thus their contribution vanishes 
in the commutative limit. 
Thus, for simplicity,  we continue here working with the full 
Hilbert space $\CH_N$, and we come back to this point 
in the discussion.

In this way, we obtain Connes' triple $(\CA_N, \BD, \CH_N)$.
We thus apply the construction of the differential calculus 
\cite{ConnesNCG}. See also \cite{ChamFroh93}.

The exterior derivative can be
defined for any element $\Ba \in \CA_N$ as:
\be
\pi(d\Ba)=[\BD,\pi(\Ba)]
\ee
where $\pi$ is the representation of the algebra $\CA_N$ in $\CH_N$. In the
following we do not distinguish between algebra element and its representation 
as far as it is obvious.
Note that in our convention $(d\Ba)^*=-d\Ba^*$.
We define then the space of $1$-forms $\Omega^1$ by
\be
\Omega^1=\{\omega |\ \omega
=\sum_i\Ba_i[\BD,\Bb_i]\ ;\ \Ba_i,\Bb_i\in \CA_N\} \ .
\ee
Thus, in the above construction the exterior derivative $d$ is a map:
\be
d\quad :\quad \CA_N\ \rightarrow\  
M_{\complex}(2)\otimes(\CA_N\otimes \CA_N^{\Bo})\ .
\ee

To define the $p$-forms, one first introduces the 
universal differential algebra 
$\widetilde{\Omega}^*=\oplus \widetilde{\Omega}^p$,
where a $p$-form is 
defined by the product
\be
\widetilde{\omega_p}\equiv \Ba_0 [\BD,\Ba_1] [\BD,\Ba_2]
\cdots [\BD,\Ba_p] \ .
\ee
The exterior derivative is
\be
d\widetilde{\omega_p}\equiv [\BD,\Ba_0]
[\BD,\Ba_1] [\BD,\Ba_2]
\cdots [\BD,\Ba_p] \ .
\ee

To obtain the graded differential algebra $\Omega^*_D$, 
with Dirac operator $\BD$ we have to 
divide out the differential ideal $\BJ+d\BJ$, $\BJ=ker(\pi)$, of 
the representation of the algebra \cite{ConnesNCG} in the Hilbert space 
$\CH_N$, 
\be
\Omega^*_D=\widetilde{\Omega}^*/(\BJ+d\BJ) \ .
\ee
The structure of this differential calculus with Dirac operator given in 
eq.(8) is discussed further in ref.\cite{CW97a}.

\section{Scalar Field}

Let us apply this construction to define a scalar field on the fuzzy sphere.
We denote the complex scalar field 
by $\Phi\in\CA_N$, then its derivative is given by
\be
d\Phi=[\BD,\Phi]={i\over\ell}
\gamma_\chi^{\Bo}\epsilon_{ijk}\sigma_i\Bx_j^{\Bo}(\BL_k\Phi) \ .
\ee
A natural choice for the action is
\be
\BS={1 \over 2(N+1)^2}\tr_{\CH}\{(d\Phi)^* d\Phi\} \ .
\ee
The above action can be evaluated by using the formula
in eq.(\ref{trace}),
\be
\BS
={-2 \over 3 \alpha^2(N+1)}\tr_{\CF}\{[\Bx_i,\Phi^\dagger][\Bx_i,\Phi]\}\ ,
\label{Action}
\ee

Recalling that in the commutative limit,
the trace ${1\over N+1}\tr_{\CF}$ corresponds to the integration 
over the sphere: 
\be
{1\over N+1}\tr_{\CF}\sim\int{d\sigma^2\over 4\pi\ell^2} \ ,
\ee
the Langrangian is given by the 'integrand' of 
eq.(\ref{Action}). Using the angular momentum it can be written as
\be
\CL={2\over3}\sum_i |\BL_i\Phi|^2\ ,
\ee

In the limit $N\rightarrow\infty$,
the angular momentum is replaced by the Killing vector, and 
we obtain the standard Lagrangian of the scalar on the sphere.

From the equation of motion of the above scalar field, we obtain the 
Casimir operator as an Laplacian. Thus, the spectrum of the 
scalar boson on the noncommutative sphere 
coincides with the spectrum of the classical scalar boson, 
until the value where the spectrum of the scalar on the 
noncommutative sphere is truncated
at its maximum angular momentum $l_{max}=N$, as expected 
\cite{MadoreCQG92}.

\section{Discussion}

In this paper we have constructed the spectral triple using 
the new chirality operator $\gamma_\chi^{\Bo}$ on the fuzzy sphere.
The differential calculus is formulated by applying Connes' construction.
Since the algebra of the fuzzy sphere is represented by the operator
in the Hilbert space $\CF_N$, the simplest way to 
introduce the grading, for example, 
would have been to double the Hilbert space as $\CF_N\oplus\CF_N$.
However, in the present approach, we 
first define 
the fermions as an $\CA_N$-bimodule and then define the Hilbert space
$\CH_N$ of these fermions \footnote{Such a strategy is also taken in 
refs.\cite{GrosKlim96,GrosKlim96a}.}.
Then the algebra $\CA_N$ is embedded into
 $\CA_N\otimes\CA_N^{\Bo}$ the elements of which are 
operators on $\CH_N$.
 The Dirac operator as well as the chirality operator are thus constructed 
in terms of elements of 
$\CA_N\otimes \CA_N^{\Bo}$ and the 
exterior derivative is a map 
$d:\CA_N\rightarrow M_{\complex}(2)\otimes\CA_N\otimes\CA_N^{\Bo}$. 
The advantage of this approach is that we have a clear correspondence
 to the commutative case.  

On the other hand, we have seen in section 3 that the  Dirac operator 
possesses zeromodes in its spectrum which 
correspond to the state of maximal angular momentum and 
they have no classical analogue. 
In the context treated here, i.e. the scalar field on the sphere, 
the contribution of these zeromodes to the action is of order 
$1/N$ and 
therefore we obtain the standard classical theory in the commutative 
limit.   

However, it is also possible to project out the zeromodes 
as follows: 
The spin $j=N+\half$ state, 
$\Psi_{N+\half,m}$, can be constructed easily and the projection operator 
to the spin $N+\half$ state is simply given by 
\be
\CP_{N+\half}=\sum_{k=0}^{2N+1}|\Psi_{N+\half, N+\half-k}>
<\Psi_{N+\half, N+\half-k}| \ .
\ee
By using this projection operator we may obtain the Hilbert space 
$\CH'_N$ which has no zeromodes. However, in such a case we must 
reconsider the whole operator algebra under this projection, 
since the algebra $\CA_N$ is not represented in the Hilbert space 
$\CH'_N$. 
Such a projection may also be avoidable by considering 
the supersymmetry \cite{GrosKlim96}.

The algebra treated here is simply a matrix algebra and by construction,
we obtain the algebra of the spherical harmonics in the limit 
$N\rightarrow \infty$, which is the reason for calling it 
the algebra of the fuzzy sphere. 
As an algebra, $\CA_N$ 
is equivalent to the algebra $M_{\complex}(N+1)$ of the 
nonsingular $(N+1)\times (N+1)$ matrices 
for given finite $N$. 

The same algebra can be obtained by applying 
Berezin's quantization for the Poisson algebra on the sphere.
By the same procedure, we can also define the fuzzy torus and
fuzzy Riemann surfaces in general, since Berezin's quantization procedure
can be applied to any K\"ahler manifold
\cite{Coburn92,KlimLes92,RawnCahen90,BordMein94,Schlich96}.  
As a result we obtain 
a noncommutative algebra which is the quantization of the function
algebra over the corresponding manifold (considering the complex manifold as a
Poisson manifold, i.e. as a phase space).  For example the
algebra of the fuzzy torus (or the noncommutative torus) 
is generated by two element
$\BS, \BT$ with the commutation relation
\be
\BS\BT=q\BT\BS
\ee 
For each case, we obtain a matrix algebra. 

Therefore, we encounter the question: if we are given a matrix algebra,
how do we find the corresponding classical geometry.  To our knowledge,
there are no general criteria to answer this question. 
It depends on the way how the limit is taken.  

When we define the Dirac operator in this paper and also in
ref.\cite{CW96}, we have required 
that the commutative limit of this operator coincides with the
standard Dirac operator on the sphere. 
However, as we disscussed above,
we also have to specify the condition of how the limit has to be taken.
One of the criteria 
to specify the commutative limit, including the differential 
operator, is the boundedness of the 
operator $[\BD, \Ba]$ for $\Ba\in\CA_N$.  
In Connes' construction, 
the boundedness of  the commutator $[\BD, \Ba]$
is required.  
It is also clear that if this commutator is not bounded, then the limit 
operator is not an 
operator in the algebra of the spherical harmonics.  
So this requirement is necessary to specify the limit as a sphere.
Furthermore, the condition of boundedness 
should be also sufficient to specify the
topology of the Riemann surface obtained in the limit 
$N\rightarrow \infty$ of the matrix algebra $M_{\complex}(N)$,
since in the present case the boundedness of the commutator
$[\BD, \Ba]$ defines a finite Poisson bracket of the algebra elements with 
the generators of the algebra $\Bx_i$.

As we have seen, 
the definition of the Dirac operator has an ambiguity 
due to the operator ordering.  In this paper we discussed only
the fuzzy sphere. However, since it is possible to discuss any 
fuzzy Riemann surface with the Berezin-Toepliz quantization of the
corresponding Riemann surface\cite{Coburn92,KlimLes92,RawnCahen90,BordMein94}, 
we also want to have a 
common description of the differential calculus on 
fuzzy Riemann surfaces which include the case of higher genus.
Such a requirement may restrict the Dirac operator. 
Another interesting aspect to consider 
all Riemann surfaces within a common framework is that
a field theory on the fuzzy Riemann surface
may be considered as a fuzzy string, the world sheet of which is
discribed by a fuzzy Riemann surface.

\bigskip
\noindent{\Large \bf{Acknowledgement}}

The authors would like to thank H. Ishikawa for helpful discussions. U.C. 
would like to acknowledge the Japan Society for Promotion of Science 
for support during the main part of this investigation.
This work is also supported by the Grant-in-Aid of 
Monbusho (the Japanese Ministry of Education, Science, Sports and Culture) 
\#09640331.


\medskip
\begin{thebibliography}{99}

\bibitem{BFSS}{T. Banks, W. Fischler, SH. Shenker, L.Susskind,
"M theory as a matrix model", Nucl.Phs.{\bf B497} 41--55 (1997).}


\bibitem{IKKT}{N. Ishibashi, H. Kawai, Y. Kitazawa, A. Tsuchiya, 
"A large-N reduced model as superstring", 
Nucl. Phys. {\bf B498} 467--491 (1997).}


\bibitem{ConnesNCG}{ A. Connes, "Noncommutative Geometry", Academic Press 
(1994).}

\bibitem{Haldane83}{F.D.M. Haldane, "Fractional Quantization of the Hall 
Effect: A Hierachy of Incompressible Quantum Fluid States", Phys. Rev. Lett. 
{\bf 51}, 605 (1983).}

\bibitem{FanoOrto86}{G. Fano, F. Ortolani, E. Colombo,
"Configuration-interaction calculations on the fractional quantum Hall 
effect", Phys. Rev. {\bf B34}, 2670--2680 (1986).}



\bibitem{Hoppe89}{J. Hoppe, "Quantum Theory of a Massless Relativistic 
Surface and a Two-Dimensional Bound State Problem", PhD Thesis, MIT (1982) 
published in 
Soryushiron Kenkyu (Kyoto) Vol.{\bf80}, 145--202 (1989).}


\bibitem{DeWHoNi88}{B. de Wit, J. Hoppe, H. Nicolai, "On the Quantum 
Mechanics of Supermembranes", Nucl. Phys. {\bf B305}, 545--581 (1988).}


\bibitem{MadoreCQG92}{J. Madore, "The fuzzy sphere", 
Class. Quant. Grav. {\bf 9}, 69-87 (1992).}


\bibitem{Berezin74}
{F.A. Berezin, "Quantization", Math. USSR Izvestija {\bf 8}, 
1109--1165 (1974).}

\bibitem{Berezin75a}{F.A. Berezin, "General concept of quantization", 
Commun. Math. Phys. {\bf 40}, 153--174 (1975).}


\bibitem{DuboKern90}{M. Dubois-Violette, R. Kerner, J. Madore,
"Noncommutative differential geometry of matrix algebras", J. Math. Phys. 
{\bf 31}, 316 (1990).}

\bibitem{GrosMado92}{H. Grosse, J. Madore, "A noncommutative 
version of the Schwinger model", Phys. Lett. {\bf 283}, 218--222 (1992).}


\bibitem{GrosKlim96}{H. Grosse, C. Klimcik, P. Presnajder, 
"Towards a finite Quantum Field Theory in Noncommutative Geometry", 
Int. Jour. Theor. Phys. {\bf 35}, 231--244 (1996).
}


\bibitem{GrosPres95}{H. Grosse, P. Presnajder, "The Dirac Operator 
on the Fuzzy Sphere", Lett. Math. Phys. {\bf 33}, 171--181 (1995).}



\bibitem{GrosKlim96a}{H. Grosse, C. Klimcik, P. Presnajder, 
"Topological Nontrivial Field Configurations in Noncommutative Geometry",
 Commun. Math. Phys. {\bf 178}, 507--526 (1996).
}


\bibitem{CW96}{U. Carow-Watamura, S. Watamura, 
           "Dirac and Chirality operator on noncommutative sphere",
          Commun. Math. Phys. {\bf 183}, 365--382 (1997).}

\bibitem{DuboMass95}{M. Dubois-Violette, T. Masson, "On the first Order 
Operators in Bimodules", Lett.Math.Phys.{\bf 37}, 467--474 (1996).}


\bibitem{CW97a}{U. Carow-Watamura, S. Watamura, in preparation.}

\bibitem{ChamFroh93}{A.H. Chamseddine, J. Fr\"ohlich, "Some Elements 
of Connes' Noncommutative Geometry, and Spacetime Geometry",
Preprint ETH-TH-93-24 (93,rec.Jul.), hep-th/9307012.} 


\bibitem{Coburn92}{L.A. Coburn, "Deformation Estimates for 
the Berezin-Toeplitz Quantization", Commun. Math. Phys. 
{\bf 149}, 415--424 (1992).}

\bibitem{KlimLes92}{S. Klimek, A. Lesniewski, 
"Quantum Riemann Surfaces I. The Unit Disk", 
Commun. Math. Phys. {\bf 146}, 103 (1992).}

\bibitem{RawnCahen90}{M. Cahen, S. Gutt, J. Rawnsley, "Quantization of 
K\"ahler Manifolds II", Transactions of the American Math. Soc. {\bf 337}, 
73--98, (1993).}


\bibitem{BordMein94}{M. Bordemann, E. Meinrenken, M. Schlichenmaier,
"Toeplitz Quantization of K\"ahler Manifolds and $gl(N),N\rightarrow\infty$ 
Limits", Commun. Math. Phys. {\bf 165}, 281--296 (1994).} 

\bibitem{Schlich96}{M. Schlichenmaier,"Berezin-Toeplitz Quantization 
of Compact K\"ahler Manifolds", q-alg/9601016.}



\end{thebibliography}
\end{document}